\documentclass[11pt, letterpaper]{article}
\usepackage[cp850]{inputenc}
\usepackage[T1]{fontenc}
\usepackage[english]{babel}
\usepackage{epsfig,graphicx,color}
\usepackage[f]{esvect}
\usepackage{amsmath, amsthm, amssymb}
\usepackage{mathrsfs}
\usepackage{url}
\usepackage{commath}
\usepackage{cite}
\usepackage{pgfplots}
\usepackage{tikz}
\usepackage{rotfloat}
\usepackage{rotating}
\usepackage{threeparttable}
\usepackage{xfrac}
\usepackage[tableposition=top]{caption}
\usepackage[margin = 1in]{geometry}

\newtheorem{theorem}{Theorem}

\allowdisplaybreaks
\def\tablenotes{\bgroup\parfillskip=0pt plus 1fil
\leftskip=0pt\relax \rightskip=0pt
\vskip2pt\footnotesize}
\def\endtablenotes{\vskip1pt\egroup}

\begin{document}
\title{Double Exponential Transformation For Computing Three-Centre Nuclear Attraction Integrals}

\author{Jordan Lovrod and Hassan Safouhi\footnote{Corresponding author: hsafouhi@ualberta.ca. \newline The corresponding author acknowledges the financial support for this research by the Natural Sciences and Engineering Research Council of Canada~(NSERC) - Grant RGPIN-2016-04317.}\\
Mathematical Division\\
Campus Saint-Jean, University of Alberta\\
8406  91 Street, Edmonton (AB) T6C 4G9, Canada}

\date{}

\maketitle
{\bf \large Abstract}.

Three-centre nuclear attraction integrals, which arise in density functional and \textit{ab initio} calculations, are one of the most time-consuming computations involved in molecular electronic structure calculations. Even for relatively small systems, millions of these laborious calculations need to be executed. Highly efficient and accurate methods for evaluating molecular integrals are therefore all the more vital in order to perform the calculations necessary for large systems. When using a basis set of $B$ functions, an analytical expression for the three-centre nuclear attraction integrals can be derived via the Fourier transform method. However, due to the presence of the highly oscillatory semi-infinite spherical Bessel integral, the analytical expression still remains problematic.  By applying the $S$ transformation, the spherical Bessel integral can be converted into a much more favorable sine integral. In the present work, we then apply two types of double exponential transformations to the resulting sine integral, which leads to a highly efficient and accurate quadrature formulae.  This method facilitates the approximation of the molecular integrals to a high predetermined accuracy, while still keeping the calculation times low. The fast convergence properties analyzed in the numerical section illustrate the advantages of the method.

\vspace*{0.75cm}

\clearpage

\section{Introduction}
The Schr\"odinger equation, first published in Schr\"odinger's 1926 paper, was a ground-breaking development that provided new insight into the behavior of the electrons in a given system \cite{Schrodinger}. Solutions to the Schr\"odinger equation yield valuable details regarding the electron placement, total energy, and other properties of the system. Nevertheless, \textit{ab initio} calculations, that is to say calculations that use the positions of the nuclei and the number of electrons to solve the Schr\"odinger equation, still present significant computational difficulties. Although the Schr\"odinger equation can be solved analytically for hydrogen and other atomic nuclei with only one electron, the inter-electron interaction involved in systems comprised of more complex atoms makes the process of analytically solving the Schr\"odinger equation extremely laborious, if not impossible.

Three-centre nuclear attraction integrals are the rate determining step of \textit{ab initio} and density functional theory (DFT) molecular structure calculations. Three-centre integrals, one of the most common types of molecular multi-centre integrals, arise when each of the three nuclei occupy a different position in space.

A common approach for calculating molecular orbitals (MOs) is to build each MO as a linear combination of atomic orbitals (LCAO), where atomic orbitals (AOs) are the solutions to the Schrodinger equation for the hydrogen atom or a hydrogen-like ion. This technique is often referred to as the LCAO-MO method. The complexity of the three-centre nuclear attraction integrals, as well as that of other multi-centre integrals, can be mitigated by strategically choosing the basis of atomic orbitals with which to perform the calculations. It has been shown that a suitable atomic orbital basis should satisfy two conditions for analytical solutions to the Schr\"odinger equation: exponential decay at infinity \cite{Agmon-85} and a cusp at the origin \cite{Kato-10-151-57}.

One choice for the basis functions are the exponential type functions (ETFs) which behave like $\exp(-x)$. ETFs are appropriate wave functions in that they satisfy the aforementioned conditions for analytical solutions to the Schr\"odinger equation. Generally, ETFs are a better suited basis than Gaussian type functions (GTFs) \cite{Boys-200-542-50, Boys-201-125-50}, which are AOs that behave like $\exp(-x^2)$. That being said however, GTFs are commonly used for chemical calculations \cite{Boys-258-402-60}. Unfortunately, these functions display neither exponential decay at infinity nor a cusp at the origin. In order to compensate, a relatively large basis set of GTFs is required. ETFs are more favorable than GTFs therefore, since a smaller basis set can be used to obtain results that are just as accurate \cite{Safouhi39}.

The most common ETFs are the Slater type functions (STFs) \cite{Slater-42-33-32}. Despite the relative simplicity of their analytical expression, however, the fact that multi-centre integrals over STFs can be extremely problematic for polyatomic molecules has averted the use of STFs as a basis set.

$B$ functions are another class of ETFs, and their use for this particular purpose was proposed by Shavitt due to their remarkably simple Weierstrass transform \cite{Shavitt-63}. $B$ functions involve reduced Bessel functions and, in fact, they can be expressed as linear combinations of STFs \cite{Filter-Steinborn-18-1-78}. Compared to other ETFs, $B$ functions are much better adapted to the evaluation of multi-centre integrals \cite{Filter-Steinborn-18-1-78, Filter-Steinborn-19-79-78, Weniger-Steinborn-28-2026-83, Weniger-1125-70-05}. This is in part due to their exceptionally straightforward Fourier transform \cite{Weniger-Steinborn-78-6121-83, Niukkanen-25-941-84}.

The Fourier transformation method introduced by Bonham \textit{et al} in 1964 \cite{Bonham-Peacher-Cox-40-3083-64}, generalized by Trivedi \textit{et al} in 1983 \cite{Trivedi-Steinborn-27-670-83}, and further generalized by Grotendorst \textit{et al} in 1988 \cite{Grotendorst-Steinborn-38-3875-88}, is one of the most successful approaches yet for the evaluation of multi-centre integrals. This particular method is especially suitable for a basis of $B$ functions. In particular, the Fourier transformation method allows for the development of an analytical expression for three-centre nuclear attraction integrals.

The challenge, however, is that these analytical expressions involve semi-infinite spherical Bessel integrals, which are highly oscillatory due to the presence of the spherical Bessel function $j_{\lambda}(vx)$. Approximating oscillatory integrals can be problematic, especially when the oscillatory part is a spherical Bessel function rather than a simple trigonometric function. Moreover, as $\lambda$ and $v$ become large, the zeros of the integrand become closer and the oscillations become stronger, thus further complicating the numerical evaluation of the integral. It is possible to rewrite such integrals as slowly convergent infinite series of integrals of alternating sign. It may then seem appropriate to apply series convergence acceleration techniques, such as Wynn's Epsilon Method \cite{Wynn-10-91-56} or Levin's $u$ transform \cite{Levin-B3-371-73}. In the case where $\lambda$ and $v$ are large, however, the excessive calculation times prevent this approach from rendering accurate results.

As shown in previous work by Safouhi~\cite{Safouhi9}, through a reformalized integration by parts with respect to $x \mathrm{d}x$, the semi-infinite spherical Bessel integral involved in the analytical expression of the three-centre nuclear attraction integrals can be transformed into an integral involving the simple sine function. This transformation, called the $S$ transformation, results in a sine integral that has equidistant zeros, thus making it much more numerically favorable than the initial spherical Bessel integral~\cite{Safouhi9, Safouhi18}.

The $S$ transformation supplies remarkable theoretical and computational power in computing spherical Bessel integrals. It applies considerable pressure on the envelope of oscillatory integrals in the asymptotic limit. While this method is highly accurate and efficient, there are some ranges of parameters where either failure is inevitable or the computation becomes extremely heavy.

In this work, after utilizing the Fourier transformation method to obtain an analytical expression, followed by the $S$ transformation to simplify the integrand to a sine function, a double exponential transformation (or DE transformation) is applied. The DE transformations introduced by Ooura \textit{et al} in 1991 \cite{Ooura-Mori-38-353-91} and refined by Ooura \textit{et al} in 1999 \cite{Ooura-Mori-112-229-99}, map the interval $(0,\infty)$ onto $(-\infty,\infty)$ and are such that the transformed integrand decays double exponentially at both infinities. The trapezoidal formula with an equal mesh size is then applied to the integral, and the resulting summation can be truncated at some moderate positive and negative values. The final result is a highly efficient quadrature formula in which relatively few function evaluations are required in order to obtain highly accurate approximations. Both the original and the refined DE transformations are applied to the spherical Bessel integral involved in the three-centre molecular integral, and their efficiency for this particular problem is compared in the discussion section. The numerical results illustrate the high accuracy of these algorithms applied to three-center nuclear attraction integrals over $B$ functions with a miscellany of different parameters.

\section{General definitions and properties}
The spherical Bessel function $j_\lambda(x)$ of order $\lambda \in \mathbb{N}_0$ is defined by \cite{Arfken-Weber-95, Abramowitz-Stegun-65}:
\begin{equation}
j_\lambda(x)  =  (-1)^\lambda x^\lambda \left(\frac{\mathrm{d}}{x \mathrm{d}x}\right)^\lambda j_0(x)
=  (-1)^\lambda x^\lambda \left(\frac{\mathrm{d}}{x \mathrm{d}x}\right)^\lambda\left(\frac{\sin{x}}{x}\right).
\end{equation}

The reduced Bessel function $\hat{k}_{n+\frac{1}{2}}(z)$ is defined by \cite{Shavitt-63, Steinborn-Filter-38-273-75}:
\begin{equation}
\hat{k}_{n+\frac{1}{2}}(z)  =  \sqrt{\frac{2}{\pi}} z^{n+\frac{1}{2}} K_{n+\frac{1}{2}}(z)
=  z^n \, e^{-z} \, \sum_{j=0}^{n} \frac{(n+j)!}{j! \,(n-j)!} \frac{1}{(2\,z)^{j}},
\end{equation}
where $n \in \mathbb{N}$, and where $K_{n+\frac{1}{2}}(z)$ is the modified Bessel function of the second kind \cite{Watson-44}.

The $B$ functions are defined by \cite{Steinborn-Filter-38-273-75, Filter-Steinborn-18-1-78}:
\begin{equation}
\label{B_functions}
B_{n,l}^{m} (\zeta, \vec{r}) = \frac{(\zeta r)^l}{2^{n+l}(n+l)!} \hat{k}_{n-\frac{1}{2}}(\zeta r) Y_l^m (\theta_{\vec{r}}, \varphi_{\vec{r}}),
\end{equation}
where $n, l, m$ are the quantum numbers and $Y_l^m (\theta, \varphi)$ denotes the surface spherical harmonic, which is defined using the Condon-Shortley phase convention \cite{Condon-Shortley-51}:
\begin{equation}
	Y_l^m (\theta, \varphi) = \textnormal{i}^{m+\abs{m}} \left[\frac{(2l+1)(l-\abs{m})!}{4\pi(l+\abs{m})!}\right] ^{1/2}
	P_l^{\abs{m}} (\cos{\theta}) \textnormal{e}^{\textnormal{i}m\varphi},
\end{equation}
where $P_l^{m}(x)$ is the associated Legendre polynomial of the $l$th degree and $m$th order.

The Fourier transform $\bar{B}_{n,l}^m(\zeta, \vec{p})$ of $B_{n,l}^{m} (\zeta, \vec{r})$ \eqref{B_functions} is given by \cite{Weniger-Steinborn-78-6121-83}
\begin{equation}
	\bar{B}_{n,l}^m = \sqrt{\frac{2}{\pi}} \zeta^{2n+l-1} \frac{(-\textnormal{i} \abs{p})^l}{(\zeta^2 + \abs{p}^2)^{n+l+1}} Y_l^m (\theta_{\vec{p}}, \varphi_{\vec{p}}).
\end{equation}

The Gaunt coefficients are defined by \cite{Gaunt-228-151-29,Homeier-Steinborn-368-31-96,Weniger-Steinborn-25-149-82,Xu-139-137-98}:
\begin{equation}
	\left\langle l_1m_1 \vert l_2m_2 \vert l_3m_3 \right\rangle = \int_{\theta=0}^{\pi} \int_{\varphi=0}^{2\pi} \left[Y_{l_1}^{m_1} (\theta, \varphi)\right]^* Y_{l_2}^{m_2} (\theta, \varphi) Y_{l_3}^{m_3} (\theta, \varphi) \sin(\theta) \mathrm{d}\theta \mathrm{d}\varphi.
\end{equation}

The Gaunt coefficients linearize the product of two spherical harmonics:
\begin{equation}
	\left[Y_{l_1}^{m_1} (\theta, \varphi)\right]^* Y_{l_2}^{m_2} (\theta, \varphi)=\sum_{l = l_{\min}, 2}^{l_1+l_2} \left\langle l_2,m_2 \vert l_1,m_1 \vert l, m_2-m_1 \right\rangle Y_{l}^{m_2-m_1} (\theta, \varphi),
\end{equation}
where the subscript $l=l_{\min, 2}$ implies that the summation index $l$ runs in steps of 2 from $l_{\min}$ to $l_1+l_2$, and where the constant $l_{\min}$ is given by:
\begin{equation}
\label{lmin}
l_{\min} =
\begin{cases}
\max(\abs{l_1-l_2}, \abs{m_2-m_1}),		& l_1+l_2+\max(\abs{l_1-l_2}, \abs{m_2-m_1}) \textnormal{ even,}\\
\max(\abs{l_1-l_2}, \abs{m_2-m_1}) + 1,	& l_1+l_2+\max(\abs{l_1-l_2}, \abs{m_2-m_1}) \textnormal{ odd.}\\
\end{cases}
\end{equation}

The three-centre nuclear attraction integral over $B$ functions is given by:
\begin{equation}
{\cal I}_{n_1,l_1,m_1}^{n_2,l_2,m_2} =
\int
\left[B_{n_1,l_1}^{m_1} \left(\zeta_{1},\vec{R}-\vv{OA}\right)\right]^{*}
\frac{1}{|\vec{R} - \vv{OC}|}
B_{n_2,l_2}^{m_2} \left(\zeta_{2},\vec{R}-\vv{OB}\right) \mathrm{d}\vec{R},
\label{EQATTRSLATER0}
\end{equation}
where $A$, $B$ and $C$ are three arbitrary points of the euclidian space ${\cal E}_3$,
while $O$ is the origin of the fixed coordinate system. $n_1$ and $n_2$ stand for the principal quantum numbers.

After performing a translation of vector $\vv{OA}$, we can write the integral ${\cal I}_{n_1,l_1,m_1}^{n_2,l_2,m_2}$ as follows:
\begin{equation}
\label{three-centre}
I_{n_1,l_1,m_1}^{n_2,l_2,m_2} (\zeta_1, \zeta_2, \vec{R}_1, \vec{R}_2)
= \int \left[B_{n_1,l_1}^{m_1}(\zeta_1, \vec{r})\right]^*
\frac{1}{\abs{\vec{r}-\vec{R}_1}}
B_{n_2,l_2}^{m_2}(\zeta_2, \vec{r} - \vec{R}_2)
\mathrm{d}\vec{r},
\end{equation}
where $\vec{r}=\vv{R}-\vv{OA}$, $\vec{R}_1=\vv{AC}$ and $\vec{R}_2=\vv{AB}$.

\section{Double exponential transformation for the computation of semi-infinite integrals}

For more details on the DE transformations and their applications, we refer the readers to the work done by Ooura \textit{et al} \cite{Ooura-Mori-38-353-91, Ooura-Mori-112-229-99}.
The double exponential (DE) transformations present a particularly efficient method for computing semi-infinite integrals of the form:
\begin{equation}
\label{slowly_decaying_and_sin}
\int_{0}^{\infty} f_1(x)\sin(vx) \mathrm{d}x,
\end{equation}
where $f_1(x)$ is a slowly decaying function and $v$ is a constant.

We will denote the semi-infinite integral we wish to evaluate by:
\begin{equation}
I = \int_{0}^{\infty} f_0(x) \mathrm{d}x.
\end{equation}

Suppose that $f_0$ is a function such that for some $\theta,\lambda \in \mathbb{R}$ and for all $n \in \mathbb{N}, n \geq n_0$, where $n_0$ is some large natural number:
\begin{equation}
f_0(\lambda n +\theta) = 0.
\end{equation}

In other words, after some large $x$ value, $f_0(x)$ has infinite equidistant zeros occurring every $\lambda$, with a shift of $\theta$ with respect to the origin.
Now provided that $f_0(x)$ is of the form of the integrand in~\eqref{slowly_decaying_and_sin}, it can easily be determined that $\lambda = \pi \, v$ and $\theta = 0$.

Now let $\phi(t)$ be a function such that:
\begin{equation}
\label{phi_condition_12}
\phi(t) \sim t  \quad \textrm{~as} \quad t \rightarrow \infty \qquad \textrm{and} \qquad
\phi(t) \sim 0  \quad \textrm{~as} \quad t \rightarrow -\infty.
\end{equation}

Let $M$ be a large positive constant. Applying the variable transformation $x = M\phi(t)$ to the integral $I$, we get:
\begin{eqnarray}
	\label{introduce_M}
	I & = & \int_{\phi^{-1}(0)}^{\phi^{-1}(\infty)} f_0( M\phi(t)) M \phi'(t) \mathrm{d}t \\
	  & = & M \int_{-\infty}^{\infty}f_0(M\phi(t))\phi'(t) \mathrm{d}t.
\end{eqnarray}

Now applying the trapezoidal rule with a mesh size of $h$ and a shift of $\dfrac{\theta}{M}$ with respect to the origin, we obtain the approximation:
\begin{equation}
	I \approx Mh \sum_{n = -\infty}^{\infty} f_0\left(M\phi\left(nh+\frac{\theta}{M}\right)\right)\phi'\left(nh+\frac{\theta}{M}\right).
\label{summation}
\end{equation}

Let $M$ and $h$ be such that:
\begin{equation}
	\label{mh=lambda}
	M h = \lambda,
\end{equation}
and consider $f_0\left(M\phi\left(nh+\frac{\theta}{M}\right)\right)$ for large $n$:
\begin{eqnarray}
f_0\left(M\phi\left(nh+\frac{\theta}{M}\right)\right) & \sim & f_0\left(M\left(nh+\frac{\theta}{M}\right)\right)
\nonumber\\ & = & f_0(Mnh + \theta)
\nonumber\\ & = & f_0(\lambda n + \theta)
\nonumber\\ & = & 0.
\label{zero_at_large_n}
\end{eqnarray}

Given the equation \eqref{zero_at_large_n}, we can truncate the summation in \eqref{summation} at some positive $N_+\in\mathbb{Z}$. Furthermore, since by \eqref{phi_condition_12}, $\phi(t)$ tends to zero as $t \rightarrow -\infty$, it follows that we can also truncate the summation in \eqref{summation} at some negative $N_-\in\mathbb{Z}$.

As a result, we obtain the approximation:
\begin{equation}
\label{final_summation}
	I \approx Mh \sum_{n = N_-}^{N_+} f_0\left(M\phi\left(nh+\frac{\theta}{M}\right)\right)\phi'\left(nh+\frac{\theta}{M}\right),
\end{equation}
for some $N_-, N_+ \in \mathbb{Z}$.

\subsection{A DE transformation}
One transformation that satisfies the conditions in \eqref{phi_condition_12} is given by~\cite{Ooura-Mori-38-353-91}:
\begin{equation}
\label{phi_1}
\phi_1(t) = \frac{t}{1-\exp(-K \sinh{t})},
\end{equation}\
where $K$ is some positive constant.

\begin{figure}[H]
	\centering
	\begin{tikzpicture}[scale=0.750]
	\begin{axis}[
		axis lines = left,
		xlabel = $t$,
		ylabel = {$\phi_1(t)$},
		ymajorgrids=true,
		grid style=dashed,
		extra y ticks = {0.16666666666},
		extra y tick style={grid=none},
		extra y tick labels = {$\frac{1}{K}$}
		]
	\addplot [domain=0.03:2, samples=400, color=black,]
	{x/(1-exp(- 6*sinh(x)))};

	\addplot [domain=-1:-0.03, samples=200,	color=black,]
	{x/(1-exp(- 6*sinh(x)))};

	\addplot[mark = o] coordinates {(0,0.16666666666)};
	\addplot[domain=-1:-0.03,dashed,mark=none,samples=200] {0.16666666666};
	\draw [densely dashed] (0,0) -- (0,0.16666666666-.03);
	\end{axis}

	\end{tikzpicture}
	\caption{The graph of $\phi_1(t)$ given by \eqref{phi_1}, where $K=6$.}
\end{figure}
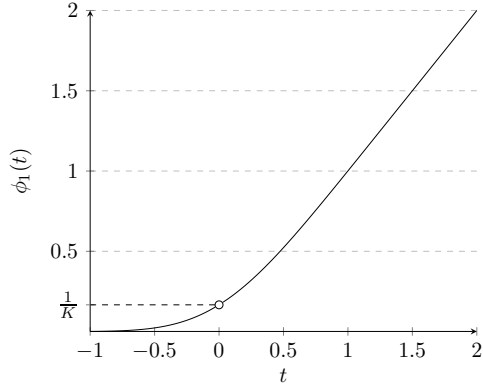

Notice that $\phi_1(t)$ encounters a singularity when $t = 0$. By Taylor expansion, the behaviour of $\phi_1(t)$ about zero can be determined. More precisely:
\begin{equation}
\phi_1(t) \sim
\left\{
\begin{array}{lll}
		0            & \quad \textrm{as} \quad & t \rightarrow -\infty\\[0.15cm]
		\dfrac{1}{K} & \quad \textrm{as} \quad &  t \rightarrow  0\\[0.15cm]
		t            & \quad \textrm{as} \quad & t \rightarrow \infty.
	\end{array}
\right.
\end{equation}

Similarly, $\phi_1'(t)$ also encounters a singularity when $t = 0$. By the same method, it can be easily found that:
\begin{equation}
\phi_1'(t) \sim
\left\{
\begin{array}{lll}
		0            & \quad \textrm{as} \quad & t\rightarrow -\infty\\[0.15cm]
		\dfrac{1}{2} & \quad \textrm{as} \quad & t\rightarrow 0\\[0.15cm]
		1            & \quad \textrm{as} \quad & t\rightarrow \infty.
	\end{array}
\right.
\end{equation}

\subsection{A more robust DE transformation}
A second transformation that satisfies the conditions in \eqref{phi_condition_12} is given by~\cite{Ooura-Mori-112-229-99}:
\begin{equation}
\label{phi_2}
\phi_2(t) = \frac{t}{1-\exp(-2x-\alpha(1-\textnormal{e}^{-t})-\beta (\textnormal{e}^t - 1))}
\end{equation}
where $\alpha, \beta$ are constants such that:
\begin{equation}
	\label{alpha-beta_conditions}
	\beta = O(1), \quad \alpha = O((M \log{M})^{-1/2}) \quad \textrm{and} \quad 0 \leq \alpha \leq \beta \leq 1.
\end{equation}

Some assignments that satisfy the constraints in \eqref{alpha-beta_conditions} are given by~\cite{Ooura-Mori-112-229-99}:
\begin{equation}
\label{alpha/beta_assignments}
\alpha = \frac{\beta}{\sqrt{1+M\log(1+M)/4\pi}} \qquad \textrm{and} \qquad  \beta = \frac{1}{4}.
\end{equation}

\begin{figure}[H]
	\centering
	\begin{tikzpicture}[scale=0.750]
	\begin{axis}[
	axis lines = left,
	xlabel = $t$,
	ylabel = {$\phi_2(t)$},
	ytick={0, 1, 1.5, 2, 2.5, 3},
	ymajorgrids=true,
	grid style=dashed,
	extra y ticks = {0.404487167214975, .5},
	extra y tick style={grid=none},
	extra y tick labels = {$\frac{1}{\alpha + \beta + 2}$, }
	]
	\draw [dashed, color = lightgray] (-2,0.5) -- (3, 0.5);

	\addplot [domain=0.04:3, samples=400, color=black,]
	{x/(1-exp(- 2*x - (0.25 / sqrt(1+ 4.5*(ln{5.5}/(4*pi))))*(1-exp(-x))-(1/4)*(exp(x)-1)))};

	\addplot [domain=-2:-0.04, samples=400, color=black,]
	{x/(1-exp(- 2*x - (0.25 / sqrt(1+ 4.5*(ln{5.5}/(4*pi))))*(1-exp(-x))-(1/4)*(exp(x)-1)))};

	\addplot[mark = o]
	coordinates {(0 , 1/2.472266319066)};

	\addplot[domain=-2:-0.04,dashed,mark=none,samples=200] {1/2.472266319066};
	\draw [densely dashed] (0,0) -- (0, 1/2.472266319066 -.04);
	\end{axis}

	\end{tikzpicture}
	\caption{The graph of $\phi_2(t)$ given by \eqref{phi_2}, with the assignments for $\alpha$ and $\beta$ in \eqref{alpha/beta_assignments}, where $M = 4.5$.}
\end{figure}
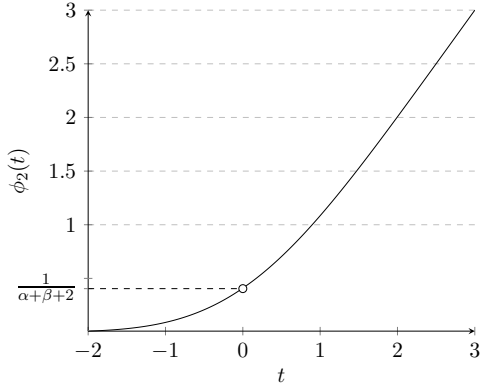

In a similar fashion to $\phi_1(t)$, $\phi_2(t)$ and $\phi_2'(t)$ both have singularities when $t = 0$. By Taylor expansion, it can be determined that:
\begin{equation}
\phi_2(t) \sim
\left\{
\begin{array}{lll}
		0                         & \quad \textrm{as} \quad & t\rightarrow -\infty\\[0.15cm]
		\dfrac{1}{2+\alpha+\beta} & \quad \textrm{as} \quad &  t\rightarrow 0\\[0.15cm]
		t                         & \quad \textrm{as} \quad & t\rightarrow \infty,\\
	\end{array}
\right.
\end{equation}
and that:
\begin{equation}
\phi_2'(t) \sim
\left\{
\begin{array}{lll}
0                                                              & \quad \textrm{as} \quad & t\rightarrow -\infty\\[0.15cm]
\dfrac{(\alpha+\beta+2)^2+(\alpha-\beta)}{2(\alpha+\beta+2)^2} & \quad \textrm{as} \quad & t\rightarrow 0\\[0.15cm]
1                                                              & \quad \textrm{as} \quad & t\rightarrow \infty.\\
\end{array}
\right.
\end{equation}

\section{The three-centre nuclear attraction integrals}
Using the Fourier transform method, an analytical expression for three-centre nuclear attraction integrals over $B$ functions \eqref{three-centre} can be found \cite{Trivedi-Steinborn-27-670-83, Grotendorst-Steinborn-38-3875-88}:
\begin{eqnarray}
\mathcal{I}_{n_1, l_1, m_1}^{n_2, l_2, m_2} & = & \frac{8(4\pi)^2(-1)^{l_1+l_2}(2l_1+1)!!(2l_2+1)!!(n_1+l_1+n_2+l_2+1)!\zeta_1^{2n_1+l_1-1}\zeta_2^{2n_2+l_2-1}}{(n_1+l_1)!(n_2+l_2)!} \nonumber\\
 & \times  & \sum_{l'_1 = 0}^{l_1} \sum_{m'_1 = -l'_1}^{l_1'} \textnormal{i}^{l_1+l'_1} (-1)^{l'_1} \frac{\left\langle l_1m_1 \vert l'_1m'_1 \vert l_1 - l'_1m_1 - m'_1\right\rangle }{(2l'_1+1)!![2(l_1-l'_1)+1]!!} \nonumber\\
 & \times  & \sum_{l'_2 = 0}^{l_2} \sum_{m'_2 = -l'_2}^{l_2'} \textnormal{i}^{l_2+l'_2} (-1)^{l'_2} \frac{\left\langle l_2m_2 \vert l'_2m'_2 \vert l_2 - l'_2m_2 - m'_2\right\rangle }{(2l'_2+1)!![2(l_2-l'_2)+1]!!} \nonumber\\
 & \times  & \sum_{l = l'_{\min}, 2}^{l'_2+l'_1} \left\langle  l'_2 m'_2 \vert l'_1 m'_1 \vert l m'_2 - m'_1 \right\rangle R_2^l Y_l^{m'_2-m'_1}(\theta_{\vec{R}_2}, \varphi_{\vec{R}_2}) \nonumber\\
 & \times & \sum_{\lambda = l''_{\min}, 2}^{l_2-l'_2+l_1-l'_1} (-\textnormal{i})^\lambda \left\langle l_2-l'_2m_2-m'_2 \vert l_1 - l'_1m_1 - m'_1 \vert \lambda\mu \right\rangle \nonumber\\
 & \times & \sum_{j = 0}^{\Delta l} \binom{\Delta l}{j} \frac{(-1)^j}{2^{n_1+n_2+l_1+l_2-j+1}(n_1+n_2+l_1+l_2-j+1)!} \nonumber\\
 & \times & \int_{s = 0}^{1} s^{n_2+l_2+l_1-l'_1}(1-s)^{n_1+l_1+l_2-l'_2} Y_\lambda^\mu(\theta_{\vec{v}}, \varphi_{\vec{v}})\nonumber\\
 & \times & \left[ \int_{x = 0}^{+ \infty} x^{n_x} \frac{\hat{k}_\nu \left[R_2 \gamma(s, x)\right]}{\left[\gamma(s,x)\right]^{n_\gamma}} j_\lambda(vx) \mathrm{d}x\right] \mathrm{d}s,
 \label{three-centre_analytical}
\end{eqnarray}
where the constant $l_{\min}$ is defined in \eqref{lmin} and:
\begin{eqnarray*}
    \gamma(s, x)  &=& \sqrt{(1-s)\zeta_1^2+s\zeta_2^2 + s(1-s)x^2}	\\
	\vec{v}		&=& (1-s)\vec{R}_2 - \vec{R}_1, \quad v \,=\, \lVert \vec{v} \rVert	  \quad \textrm{and} \quad
	 R_2\,=\, \Vert \vec{R}_2 \rVert				\\
	n_x							&=& l_1-l'_1+l_2-l'_2					\\
	n_\gamma							&=& 2(n_1+l_1+n_2+l_2) - (l'_1+l'_2) - l+1	\\
	\nu							&=& n_1+n_2+l_1+l_2-l-j+ \frac{1}{2}		\\
	\mu							&=& (m_2-m'_2) - (m_1-m'_1)				\\
	\Delta l					&=& \left[(l'_1+l'_2-l)/2\right].
\end{eqnarray*}

The semi-infinite integral in \eqref{three-centre_analytical}, which will from now on be referred to as $\mathcal{I}(s)$, is defined by:
\begin{equation}
\label{I(s)}
\mathcal{I}(s) = \int_{0}^{+ \infty} x^{n_x} \frac{\hat{k}_\nu \left[R_2 \gamma(s, x)\right]}{\left[\gamma(s,x)\right]^{n_\gamma}} j_\lambda(vx) \mathrm{d}x.
\end{equation}

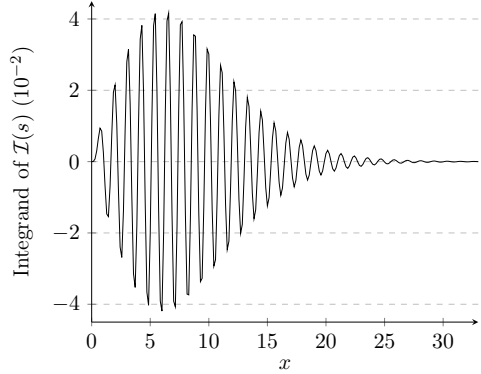
\begin{figure}[H]
	\centering
	\begin{tikzpicture}[scale=0.750]
	\begin{axis}[
	axis lines = left,
	xlabel={$x$},
	ylabel={Integrand of  $\mathcal{I}(s)$ ($10^{-2}$)},
	xmin=0, xmax=33,
	ymin=-4.5, ymax=4.5,
	ymajorgrids=true,
	grid style=dashed,
	]
	\addplot[no marks]
	coordinates {
	( 0.000, +0.0000)	( 0.143, +0.0054)	( 0.286, +0.0747)	( 0.429, +0.2972)	( 0.571, +0.6528)	( 0.714, +0.9367)	( 0.857, +0.8553)	( 1.000, +0.2543)
	( 1.143, -0.6862)	( 1.286, -1.4668)	( 1.429, -1.5396)	( 1.571, -0.6884)	( 1.714, +0.7396)	( 1.857, +1.9501)	( 2.000, +2.1477)	( 2.143, +1.0594)
	( 2.286, -0.8097)	( 2.429, -2.3954)	( 2.571, -2.6858)	( 2.714, -1.3702)	( 2.857, +0.8926)	( 3.000, +2.7962)	( 3.143, +3.1476)	( 3.286, +1.6195)
	( 3.429, -0.9833)	( 3.571, -3.1444)	( 3.714, -3.5271)	( 3.857, -1.8067)	( 4.000, +1.0769)	( 4.143, +3.4333)	( 4.286, +3.8204)	( 4.429, +1.9333)
	( 4.571, -1.1684)	( 4.714, -3.6584)	( 4.857, -4.0269)	( 5.000, -2.0027)	( 5.143, +1.2532)	( 5.286, +3.8177)	( 5.429, +4.1489)	( 5.571, +2.0200)
	( 5.714, -1.3274)	( 5.857, -3.9116)	( 6.000, -4.1916)	( 6.143, -1.9917)	( 6.286, +1.3876)	( 6.429, +3.9428)	( 6.571, +4.1622)	( 6.714, +1.9249)
	( 6.857, -1.4317)	( 7.000, -3.9161)	( 7.143, -4.0696)	( 7.286, -1.8274)	( 7.429, +1.4584)	( 7.571, +3.8375)	( 7.714, +3.9240)	( 7.857, +1.7068)
	( 8.000, -1.4672)	( 8.143, -3.7143)	( 8.286, -3.7356)	( 8.429, -1.5703)	( 8.571, +1.4586)	( 8.714, +3.5541)	( 8.857, +3.5148)	( 9.000, +1.4244)
	( 9.143, -1.4336)	( 9.286, -3.3649)	( 9.429, -3.2713)	( 9.571, -1.2749)	( 9.714, +1.3938)	( 9.857, +3.1545)	(10.000, +3.0143)	(10.143, +1.1266)
	(10.286, -1.3414)	(10.429, -2.9302)	(10.571, -2.7515)	(10.714, -0.9834)	(10.857, +1.2785)	(11.000, +2.6984)	(11.143, +2.4898)	(11.286, +0.8481)
	(11.429, -1.2074)	(11.571, -2.4651)	(11.714, -2.2346)	(11.857, -0.7228)	(12.000, +1.1305)	(12.143, +2.2350)	(12.286, +1.9902)	(12.429, +0.6087)
	(12.571, -1.0500)	(12.714, -2.0121)	(12.857, -1.7598)	(13.000, -0.5065)	(13.143, +0.9678)	(13.286, +1.7994)	(13.429, +1.5456)	(13.571, +0.4162)
	(13.714, -0.8856)	(13.857, -1.5992)	(14.000, -1.3488)	(14.143, -0.3376)	(14.286, +0.8050)	(14.429, +1.4129)	(14.571, +1.1700)	(14.714, +0.2700)
	(14.857, -0.7270)	(15.000, -1.2415)	(15.143, -1.0091)	(15.286, -0.2127)	(15.429, +0.6527)	(15.571, +1.0852)	(15.714, +0.8658)	(15.857, +0.1646)
	(16.000, -0.5827)	(16.143, -0.9440)	(16.286, -0.7390)	(16.429, -0.1248)	(16.571, +0.5174)	(16.714, +0.8174)	(16.857, +0.6277)	(17.000, +0.0923)
	(17.143, -0.4572)	(17.286, -0.7048)	(17.429, -0.5308)	(17.571, -0.0660)	(17.714, +0.4020)	(17.857, +0.6052)	(18.000, +0.4469)	(18.143, +0.0451)
	(18.286, -0.3520)	(18.429, -0.5176)	(18.571, -0.3747)	(18.714, -0.0287)	(18.857, +0.3069)	(19.000, +0.4412)	(19.143, +0.3130)	(19.286, +0.0161)
	(19.429, -0.2665)	(19.571, -0.3747)	(19.714, -0.2604)	(19.857, -0.0066)	(20.000, +0.2306)	(20.143, +0.3172)	(20.286, +0.2159)	(20.429, -0.0004)
	(20.571, -0.1989)	(20.714, -0.2678)	(20.857, -0.1784)	(21.000, +0.0054)	(21.143, +0.1709)	(21.286, +0.2254)	(21.429, +0.1469)	(21.571, -0.0088)
	(21.714, -0.1464)	(21.857, -0.1891)	(22.000, -0.1205)	(22.143, +0.0109)	(22.286, +0.1251)	(22.429, +0.1583)	(22.571, +0.0986)	(22.714, -0.0121)
	(22.857, -0.1065)	(23.000, -0.1322)	(23.143, -0.0805)	(23.286, +0.0127)	(23.429, +0.0905)	(23.571, +0.1102)	(23.714, +0.0654)	(23.857, -0.0126)
	(24.000, -0.0767)	(24.143, -0.0916)	(24.286, -0.0531)	(24.429, +0.0123)	(24.571, +0.0649)	(24.714, +0.0760)	(24.857, +0.0429)	(25.000, -0.0116)
	(25.143, -0.0547)	(25.286, -0.0629)	(25.429, -0.0346)	(25.571, +0.0108)	(25.714, +0.0461)	(25.857, +0.0520)	(26.000, +0.0278)	(26.143, -0.0100)
	(26.286, -0.0387)	(26.429, -0.0429)	(26.571, -0.0223)	(26.714, +0.0091)	(26.857, +0.0325)	(27.000, +0.0353)	(27.143, +0.0178)	(27.286, -0.0082)
	(27.429, -0.0272)	(27.571, -0.0290)	(27.714, -0.0142)	(27.857, +0.0073)	(28.000, +0.0227)	(28.143, +0.0238)	(28.286, +0.0113)	(28.429, -0.0065)
	(28.571, -0.0189)	(28.714, -0.0195)	(28.857, -0.0090)	(29.000, +0.0057)	(29.143, +0.0158)	(29.286, +0.0160)	(29.429, +0.0071)	(29.571, -0.0050)
	(29.714, -0.0131)	(29.857, -0.0130)	(30.000, -0.0056)	(30.143, +0.0044)	(30.286, +0.0109)	(30.429, +0.0106)	(30.571, +0.0044)	(30.714, -0.0038)
	(30.857, -0.0090)	(31.000, -0.0087)	(31.143, -0.0034)	(31.286, +0.0033)	(31.429, +0.0075)	(31.571, +0.0071)	(31.714, +0.0027)	(31.857, -0.0028)
	(32.000, -0.0062)	(32.143, -0.0057)	(32.286, -0.0021)	(32.429, +0.0024)	(32.571, +0.0051)	(32.714, +0.0047)	(32.857, +0.0016)	(33.000, -0.0021)
	};
	\end{axis}
	\end{tikzpicture}
	\caption{The integrand of $\mathcal{I}(s)$ where $s=0.01$, $\nu=\sfrac{9}{2}$, $n_\gamma=7$, $n_x=2$, $\lambda=2$, $R_1=9$, $\zeta_1=2$, $R_2=3.5$, and $\zeta_2=1$ (see row 9 of Table \ref{table:with_phi12}).}
	\label{Integrand_I(s)}
\end{figure}

The challenge of evaluating \eqref{three-centre_analytical} is mainly caused by the presence of the semi-infinite integral $\mathcal{I}(s)$ \eqref{I(s)}, whose integrand is highly oscillatory due to the presence of the spherical Bessel function. This is especially the case when $\lambda$ and $v$ are large. Notice also that when the value of $s$ is close to $0$ or $1$, the oscillation of the integrand become sharp. Indeed, if we make the substitution $s = 0$ or $s = 1$ in the integrand, the exponentially decaying part of the integrand:
\begin{equation*}
\frac{\hat{k}_\nu \left[R_2 \gamma(s, x)\right]}{\left[\gamma(s,x)\right]^{n_\gamma}},
\end{equation*}
becomes constant, and the integrand can be reduced to $x^{n_x}j_\lambda(vx)$. As a result, the oscillations of $j_\lambda(vx)$ cannot be restrained, and the evaluation of $\mathcal{I}(s)$ becomes most laborious.
In Figure \ref{Integrand_I(s)}, the highly oscillatory nature of the integrand can be observed.

\subsection{Application of the $S$ transformation}
We shall now reiterate a theorem which is stated and proven by Safouhi \cite{Safouhi9}.
\begin{theorem}~\cite{Safouhi9}
	\label{S-traformation_theorem}
	Suppose that $f(x)$ is a function integrable on \textnormal{[0, $\infty$[} and of the form:
	\begin{equation*}
		f(x) = g(x)j_\lambda(x),
	\end{equation*}
	where $g(x) \in C^2$ \textnormal{([0, $\infty$[)}, which is the set of all twice continuously differentiable functions defined on the interval \textnormal{[0, $\infty$[}. Now if $g(x)$ is such that $\left(\frac{d}{xdx}\right)^l (x^{\lambda-1}g(x))$ for $l=0, 1, 2, \ldots, \lambda$ are definied and:
	\begin{eqnarray}
		\lim\limits_{x \rightarrow 0} x^{l-\lambda+1} \left(\frac{d}{xdx}\right)
		\left(x^{\lambda-1} g(x)\right)
		j_{\lambda-l-1}(x) & = & 0\nonumber\\
		\lim\limits_{x \rightarrow \infty} x^{l-\lambda+1} \left(\frac{d}{xdx}\right)
		\left(x^{\lambda-1} g(x)\right)
		j_{\lambda-l-1}(x) & = & 0,
	\end{eqnarray}
	holds true for all $l=0, 1, 2, \ldots, \lambda-1$, then:
	\begin{equation*}
		\int_{0}^{\infty} f(x) \mathrm{d}x = \int_{0}^{\infty}
		\left[\left(\frac{d}{xdx}\right)^\lambda (x^{\lambda-1}g(x))\right] \sin(x) \mathrm{d}x.
	\end{equation*}
\end{theorem}

In short, this transformation, called the $S$ transformation, can simplify spherical Bessel integrals into integrals involving the sine function. As such, the integral $\mathcal{I}(s)$ can be rewritten as:
\begin{equation}
\label{simplification_tilde_I}
\mathcal{I}(s) = \frac{1}{v^{\lambda+1}} \int_{0}^{+ \infty} \left[\left( \frac{d}{xdx}\right) ^\lambda \left( x^{n_x+\lambda-1} \frac{\hat{k}_\nu \left[R_2 \gamma(s, x)\right]}{\left[\gamma(s,x)\right]^{n_\gamma}}\right)\right] \sin(vx)dx.
\end{equation}

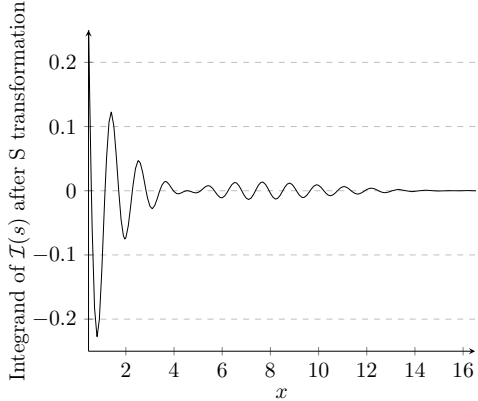
\begin{figure}[H]
	\centering
	\begin{tikzpicture}[scale=0.750]
	\begin{axis}[
	axis lines = left,
	xlabel={$x$},
	ylabel={Integrand of $\mathcal{I}(s)$ after S transformation},
	xmin=.45, xmax=16.5,
	ymin=-.25, ymax=.25,
	ymajorgrids=true,
	grid style=dashed,
	]
	\addplot[no marks]
	coordinates {
		(0.000,	+1.09507)	(0.100,	+1.03928)	(0.200,	+0.88226)	(0.300,	+0.65272)	(0.400,	+0.39139)	(0.500,	+0.14209)	(0.568,	+0.00000)	(0.600,	-0.05738)
		(0.700,	-0.18271)	(0.800,	-0.22693)	(0.808,	-0.22716)	(0.900,	-0.20015)	(1.000,	-0.12540)	(1.100,	-0.03193)	(1.135,	+0.00000)	(1.200,	+0.05222)
		(1.300,	+0.10665)	(1.389,	+0.12228)	(1.400,	+0.12206)	(1.500,	+0.10100)	(1.600,	+0.05528)	(1.700,	+0.00144)	(1.703,	+0.00000)	(1.800,	-0.04438)
		(1.900,	-0.07068)	(1.960,	-0.07498)	(2.000,	-0.07314)	(2.100,	-0.05466)	(2.200,	-0.02362)	(2.270,	+0.00000)	(2.300,	+0.00936)	(2.400,	+0.03465)
		(2.500,	+0.04627)	(2.526,	+0.04680)	(2.600,	+0.04303)	(2.700,	+0.02816)	(2.800,	+0.00776)	(2.838,	+0.00000)	(2.900,	-0.01144)	(3.000,	-0.02413)
		(3.089,	-0.02782)	(3.100,	-0.02777)	(3.200,	-0.02290)	(3.300,	-0.01249)	(3.400,	-0.00062)	(3.406,	+0.00000)	(3.500,	+0.00894)	(3.600,	+0.01385)
		(3.645,	+0.01435)	(3.700,	+0.01365)	(3.800,	+0.00959)	(3.900,	+0.00389)	(3.973,	+0.00000)	(4.000,	-0.00121)	(4.100,	-0.00420)	(4.171,	-0.00480)
		(4.200,	-0.00471)	(4.300,	-0.00339)	(4.400,	-0.00148)	(4.500,	-0.00018)	(4.541,	+0.00000)	(4.553,	+0.00001)	(4.565,	+0.00000)	(4.600,	-0.00014)
		(4.700,	-0.00122)	(4.800,	-0.00264)	(4.900,	-0.00337)	(4.905,	-0.00340)	(5.000,	-0.00262)	(5.100,	-0.00025)	(5.108,	+0.00000)	(5.200,	+0.00309)
		(5.300,	+0.00623)	(5.400,	+0.00785)	(5.420,	+0.00790)	(5.500,	+0.00705)	(5.600,	+0.00373)	(5.676,	+0.00000)	(5.700,	-0.00129)	(5.800,	-0.00650)
		(5.900,	-0.01014)	(5.972,	-0.01100)	(6.000,	-0.01088)	(6.100,	-0.00824)	(6.200,	-0.00284)	(6.243,	+0.00000)	(6.300,	+0.00376)	(6.400,	+0.00953)
		(6.500,	+0.01260)	(6.533,	+0.01280)	(6.600,	+0.01192)	(6.700,	+0.00758)	(6.800,	+0.00081)	(6.811,	+0.00000)	(6.900,	-0.00634)	(7.000,	-0.01169)
		(7.096,	-0.01360)	(7.100,	-0.01357)	(7.200,	-0.01138)	(7.300,	-0.00575)	(7.379,	+0.00000)	(7.400,	+0.00161)	(7.500,	+0.00847)	(7.600,	+0.01277)
		(7.661,	+0.01350)	(7.700,	+0.01321)	(7.800,	+0.00972)	(7.900,	+0.00337)	(7.946,	+0.00000)	(8.000,	-0.00388)	(8.100,	-0.00984)	(8.200,	-0.01277)
		(8.226,	-0.01290)	(8.300,	-0.01185)	(8.400,	-0.00744)	(8.500,	-0.00095)	(8.514,	+0.00000)	(8.600,	+0.00563)	(8.700,	+0.01035)	(8.792,	+0.01190)
		(8.800,	+0.01187)	(8.900,	+0.00984)	(9.000,	+0.00498)	(9.081,	+0.00000)	(9.100,	-0.00115)	(9.200,	-0.00671)	(9.300,	-0.01006)	(9.358,	-0.01060)
		(9.400,	-0.01032)	(9.500,	-0.00754)	(9.600,	-0.00269)	(9.649,	+0.00000)	(9.700,	+0.00273)	(9.800,	+0.00707)	(9.900,	+0.00913)	(9.924,	+0.00922)
		(10.000,+0.00842)	(10.100,+0.00528)	(10.200,+0.00078)	(10.217,+0.00000)	(10.300,-0.00369)	(10.400,-0.00683)	(10.490,-0.00780)	(10.500,-0.00778)
		(10.600,-0.00642)	(10.700,-0.00328)	(10.784,+0.00000)	(10.800,+0.00062)	(10.900,+0.00407)	(11.000,+0.00611)	(11.056,+0.00643)	(11.100,+0.00624)
		(11.200,+0.00455)	(11.300,+0.00166)	(11.352,+0.00000)	(11.400,-0.00149)	(11.500,-0.00397)	(11.600,-0.00511)	(11.622,-0.00515)	(11.700,-0.00468)
		(11.800,-0.00294)	(11.900,-0.00049)	(11.919,+0.00000)	(12.000,+0.00189)	(12.100,+0.00352)	(12.188,+0.00400)	(12.200,+0.00399)	(12.300,+0.00327)
		(12.400,+0.00167)	(12.487,+0.00000)	(12.500,-0.00025)	(12.600,-0.00191)	(12.700,-0.00286)	(12.753,-0.00299)	(12.800,-0.00289)	(12.900,-0.00209)
		(13.000,-0.00078)	(13.054,+0.00000)	(13.317,+0.00214)	(13.622,+0.00000)	(13.880,-0.00142)	(14.190,+0.00000)	(14.439,+0.00085)	(14.757,+0.00000)
		(14.985,-0.00039)	(15.325,+0.00000)	(15.464,+0.00007)	(15.648,+0.00000)	(15.777,-0.00004)	(15.892,+0.00000)	(16.223,+0.00025)	(16.460,+0.00000)
		(16.764,-0.00043)	(17.028,+0.00000)	(17.321,+0.00055)	(17.595,+0.00000)	(17.884,-0.00062)	(18.163,+0.00000)	(18.448,+0.00065)	(18.730,+0.00000)
	};
	\end{axis}
	\end{tikzpicture}
	\caption{The integrand of $\mathcal{I}(s)$ in \eqref{simplification_tilde_I} after applying the S transformation, where $s=0.01$, $\nu=\sfrac{9}{2}$, $n_\gamma=7$, $n_x=2$, $\lambda=2$, $R_1=9$, $\zeta_1=2$, $R_2=3.5$, and $\zeta_2=1$ (see row 9 of Table \ref{table:with_phi12}). The integrand prior to the transformation can be observed in Figure \ref{Integrand_I(s)}.}
	\label{Integrand_I(s)_after_S}
\end{figure}

Now, let the part of the integrand not involving the sine function be referred to as $f(x)$, i.e. let:
\begin{equation}
\label{define_f(x)}
	f(x) = \left( \frac{d}{xdx}\right) ^\lambda \left( x^{n_x+\lambda-1} \frac{\hat{k}_\nu \left[R_2 \gamma(s, x)\right]}{\left[\gamma(s,x)\right]^{n_\gamma}}\right).
\end{equation}
Notice that $f(x)$ is a decaying function and is therefore of the same form as the function $f_1(x)$, which was introduced in \eqref{slowly_decaying_and_sin}. The integrand of $\mathcal{I}(s)$ is hence a suitable candidate for a DE variable transformation, such as $\phi_1(t)$ \eqref{phi_1} or $\phi_2(t)$ \eqref{phi_2}.

Applying the DE transformation $\phi(t)$ to the equation for $\mathcal{I}(s)$  in \eqref{simplification_tilde_I}, we obtain the approximation:
\begin{eqnarray}
\mathcal{I}(s)	&    =    & \frac{1}{v^{\lambda+1}} \int_{x = 0}^{+ \infty} f(x) \sin(vx)dx \nonumber\\
						& \approx & \frac{Mh}{v^{\lambda+1}} \sum_{n = N_-}^{N_+} f\left(M\phi\left(nh+\frac{\theta}{M}\right)\right)\phi'\left(nh+\frac{\theta}{M}\right) \nonumber\\
						& = & \mathcal{I}_M(s).
						\label{approximation}
\end{eqnarray}

\section{Numerical results and discussion}

Table~\ref{table:with_phi12} lists values for $\mathcal{I}(s)$ obtained using the the approximation \eqref{approximation} with the transformations $\phi_1(t)$ \eqref{phi_1} and $\phi_2(t)$ \eqref{phi_2}. In this first table, we restrict the variable assignments to values that can be handled by a MATLAB built-in numerical integration function that uses global adaptive quadrature set to an accuracy of $15$ correct digits. In Table~\ref{table:using_phi12}, we restrict the variable assignments to values that cannot be handled by the MATLAB integration function. For more problematic variable assignments, particularly when $\lambda$ and/or $v$ are large, the MATLAB automatic integrator fails to provide results to the same accuracy efficiently. We therefore opted to use the same algorithm in Python which, with the help of the symbolic computation package SymPy, was able to complete the approximations significantly faster.
The relative errors corresponding to the approximations using each transformation are listed in Table~\ref{table:using_phi12} as $\varepsilon_M^{\phi_1}$ and $\varepsilon_M^{\phi_2}$, respectively.

As shown by Ooura \textit{et al} \cite{Ooura-Mori-38-353-91, Ooura-Mori-112-229-99}, the relative error for the approximation $\mathcal{I}_M(s)$ using a DE transformation can be written:
\begin{equation}
\label{error}
\varepsilon_M \approx \exp\left(-\frac{c}{h}\right),
\end{equation}
where $c$ is a constant. Ooura \textit{et al} argue that we can assume that the relative error is approximately given by \cite{Ooura-Mori-112-229-99}:
\begin{equation}
	\varepsilon_M \, \approx \, \exp\left(\frac{-A}{h}\right) \,=\, \exp\left(\frac{-AM}{v \, \pi}\right),
\end{equation}
where $A$ is a constant. For the transformation $\phi_1(t)$ \eqref{phi_1}, we use $A = 2$, and for the transformation $\phi_2(t)$ \eqref{phi_2}, we use $A = 5$.

\subsection{Effect of collocation points}
\begin{figure}[H]
	\centering
	\begin{tikzpicture}[scale=0.750]
	\begin{axis}[
	axis lines = left,
	xlabel={$n$ for the transformation $\phi_1(t)$},
	ylabel={Absolute error ($10^{-3}$)},
	xmin=70, xmax=120,
	ymin=-0, ymax=10,
	ymajorgrids=true,
	grid style=dashed,
	]
	\addplot[only marks, mark size=1.3pt, black, fill = white]
	coordinates {
(	70	,	9.0536631055156000	)
(	72	,	6.7534468717278000	)
(	74	,	4.3834695556151000	)
(	76	,	2.6067419485517000	)
(	78	,	1.4594331718510000	)
(	80	,	0.7813259587446000	)
(	82	,	0.4036299741182000	)
(	84	,	0.2022535532329000	)
(	86	,	0.0985796446472000	)
(	88	,	0.0467965523463000	)
(	90	,	0.0216423238195000	)
(	92	,	0.0097475071905000	)
(	94	,	0.0042717829315000	)
(	96	,	0.0018193468960000	)
(	98	,	0.0007518764540000	)
(	100	,	0.0003009694622000	)
(	102	,	0.0001164571307000	)
(	104	,	0.0000434617228000	)
(	106	,	0.0000156061203000	)
(	108	,	0.0000053777356000	)
(	110	,	0.0000017733258000	)
(	112	,	0.0000005579431000	)
(	114	,	0.0000001668497000	)
(	116	,	0.0000000473684000	)
(	118	,	0.0000000128455000	)
(	120	,	0.0000000033455000	)
(	122	,	0.0000000010591000	)
(	124	,	0.0000000000533000	)
(	126	,	0.0000000003047000	)
(	128	,	0.0000000004656000	)
(	130	,	0.0000000008434000	)
(	132	,	0.0000000005251000	)
(	134	,	0.0000000003654000	)
(	136	,	0.0000000003654000	)
(	138	,	0.0000000003654000	)
	};
	\addplot [domain=69:120, samples=200,	color=black,]
	{100*exp(-((x-47)^2)/228)};
	\end{axis}
\end{tikzpicture}~\begin{tikzpicture}[scale=0.750]
	\begin{axis}[
	axis lines = left,
	xlabel={$n$ for the transformation $\phi_2(t)$},
	xmin=45, xmax=90,
	ymin=-0, ymax=10,
	ymajorgrids=true,
	grid style=dashed,
	]
	\addplot[only marks, mark size=1.3pt, black, fill = white]
	coordinates {
	(45	,	9.1299410861758)
	(47	,	6.5887323478478)
	(49	,	4.1089337309503)
	(51	,	2.3289067750621)
	(53	,	1.2306408562994)
	(55	,	0.6142337344385)
	(57	,	0.2912613390752)
	(59	,	0.1313711904183)
	(61	,	0.0562518282602)
	(63	,	0.0227676830275)
	(65	,	0.0086559227693)
	(67	,	0.0030660276641)
	(69	,	0.0010016713957)
	(71	,	0.0002981777609)
	(73	,	0.0000797112179)
	(75	,	0.0000188132527)
	(77	,	0.0000038409399)
	(79	,	0.0000006642113)
	(81	,	0.0000000930942)
	(83	,	0.0000000116265)
	(85	,	0.0000000002727)
	(87	,	0.0000000000000)
	(97	,	0.0000000000000)
	};
	\addplot [domain=45:90, samples=200,	color=black,]
	{1000*exp(-((x-8)^2)/300)};
	\end{axis}
	\end{tikzpicture}

	\caption{Absolute error dependence of $\mathcal{I}_M(s)$ \eqref{approximation} on the total number of collocation points using transformations $\phi_1(t)$ \eqref{phi_1} and $\phi_2(t)$ \eqref{phi_2}, and the variable assignments from the second row of Table~\ref{table:with_phi12}. The absolute error was found by comparing the results of the approximation with the result obtained by a MATLAB built-in numerical integration function that uses global adaptive quadrature.}
	\label{collocation_points_phi1}
\end{figure}
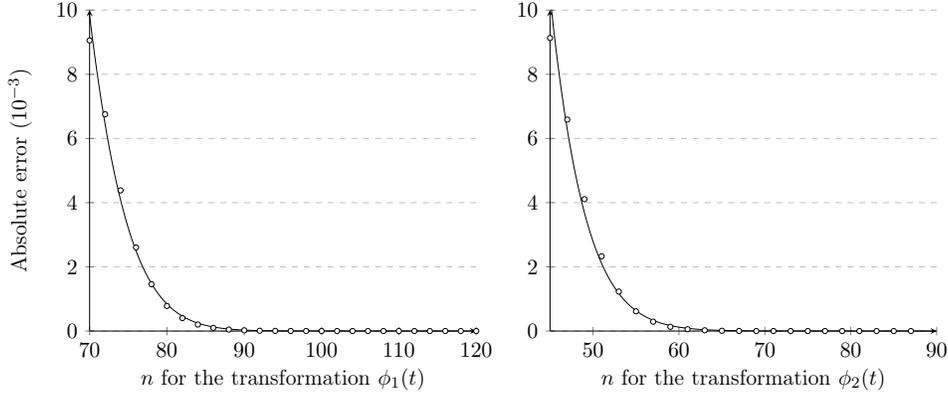

For the variable assignments used in Figure~\ref{collocation_points_phi1} for the transformation $\phi_1(t)$ \eqref{phi_1}, in our present integrator, an optimal value for $M$ was found to be $M \approx 54.25338$.  Over the course of our calculations for Figure~\ref{collocation_points_phi1}, therefore, we assigned said value to $M$ in order to preserve the behaviour of the approximation that was conducted by the integrator.

In our integrator, we first took the summation in \eqref{approximation} over positive $n$ values, until the addition of higher terms stopped significantly contributing to the overall approximation. Similarly, we then took the summation over negative $n$ values, again truncating the summation once the contributions of each term became less significant than the predetermined error tolerance. Using this method, for the variable assignments in Figure~\ref{collocation_points_phi1} and while using transformation $\phi_1(t)$, the summation in~\eqref{approximation} was truncated at -96 and at 41.

\begin{figure}[H]
	\centering
	\begin{tikzpicture}[scale=0.750]
	\begin{axis}[
	axis lines = left,
	xlabel={$x$},
	ylabel={Integrand of $\mathcal{I}(s) = \int_{-\infty}^{\infty} f\left(M\phi(t)\right) \cdot M\phi'(t) \mathrm{d} t$},
	xmin=-6, xmax=6,
	ymin=-0.8, ymax=1,
	ymajorgrids=true,
	grid style=dashed,
	]
	\addplot[no marks]
	coordinates {
		(-6.000000,  0.000000)  (-5.979967,  0.000000)  (-5.959933,  0.000000)  (-5.939900,  0.000000)  (-5.919866,  0.000000)
		(-4.898164,  0.000000)  (-4.878130,  0.000000)  (-4.858097,  0.000000)  (-4.838063,  0.000001)  (-4.818030,  0.000001)
		(-4.797997,  0.000001)  (-4.777963,  0.000001)  (-4.757930,  0.000002)  (-4.737896,  0.000002)  (-4.717863,  0.000002)
		(-4.697830,  0.000003)  (-4.677796,  0.000004)  (-4.657763,  0.000005)  (-4.637730,  0.000006)  (-4.617696,  0.000008)
		(-4.597663,  0.000010)  (-4.577629,  0.000012)  (-4.557596,  0.000015)  (-4.537563,  0.000018)  (-4.517529,  0.000022)
		(-4.497496,  0.000027)  (-4.477462,  0.000033)  (-4.457429,  0.000040)  (-4.437396,  0.000048)  (-4.417362,  0.000058)
		(-4.397329,  0.000070)  (-4.377295,  0.000084)  (-4.357262,  0.000100)  (-4.337229,  0.000118)  (-4.317195,  0.000140)
		(-4.297162,  0.000165)  (-4.277129,  0.000195)  (-4.257095,  0.000229)  (-4.237062,  0.000268)  (-4.217028,  0.000313)
		(-4.196995,  0.000364)  (-4.176962,  0.000423)  (-4.156928,  0.000490)  (-4.136895,  0.000567)  (-4.116861,  0.000653)
		(-4.096828,  0.000751)  (-4.076795,  0.000862)  (-4.056761,  0.000987)  (-4.036728,  0.001128)  (-4.016694,  0.001286)
		(-3.996661,  0.001462)  (-3.976628,  0.001660)  (-3.956594,  0.001880)  (-3.936561,  0.002125)  (-3.916528,  0.002397)
		(-3.896494,  0.002699)  (-3.876461,  0.003033)  (-3.856427,  0.003402)  (-3.836394,  0.003809)  (-3.816361,  0.004257)
		(-3.796327,  0.004749)  (-3.776294,  0.005289)  (-3.756260,  0.005880)  (-3.736227,  0.006527)  (-3.716194,  0.007234)
		(-3.696160,  0.008004)  (-3.676127,  0.008843)  (-3.656093,  0.009754)  (-3.636060,  0.010744)  (-3.616027,  0.011816)
		(-3.595993,  0.012978)  (-3.575960,  0.014233)  (-3.555927,  0.015589)  (-3.535893,  0.017051)  (-3.515860,  0.018625)
		(-3.495826,  0.020319)  (-3.475793,  0.022139)  (-3.455760,  0.024092)  (-3.435726,  0.026186)  (-3.415693,  0.028428)
		(-3.395659,  0.030827)  (-3.375626,  0.033390)  (-3.355593,  0.036126)  (-3.335559,  0.039043)  (-3.315526,  0.042151)
		(-3.295492,  0.045459)  (-3.275459,  0.048977)  (-3.255426,  0.052713)  (-3.235392,  0.056679)  (-3.215359,  0.060885)
		(-3.195326,  0.065340)  (-3.175292,  0.070057)  (-3.155259,  0.075047)  (-3.135225,  0.080319)  (-3.115192,  0.085888)
		(-3.095159,  0.091763)  (-3.075125,  0.097958)  (-3.055092,  0.104484)  (-3.035058,  0.111355)  (-3.015025,  0.118583)
		(-2.994992,  0.126181)  (-2.974958,  0.134162)  (-2.954925,  0.142540)  (-2.934891,  0.151327)  (-2.914858,  0.160537)
		(-2.894825,  0.170184)  (-2.874791,  0.180281)  (-2.854758,  0.190840)  (-2.834725,  0.201876)  (-2.814691,  0.213402)
		(-2.794658,  0.225429)  (-2.774624,  0.237970)  (-2.754591,  0.251036)  (-2.734558,  0.264639)  (-2.714524,  0.278788)
		(-2.694491,  0.293494)  (-2.674457,  0.308763)  (-2.654424,  0.324602)  (-2.634391,  0.341016)  (-2.614357,  0.358008)
		(-2.594324,  0.375579)  (-2.574290,  0.393727)  (-2.554257,  0.412447)  (-2.534224,  0.431731)  (-2.514190,  0.451566)
		(-2.494157,  0.471937)  (-2.474124,  0.492822)  (-2.454090,  0.514194)  (-2.434057,  0.536019)  (-2.414023,  0.558257)
		(-2.393990,  0.580860)  (-2.373957,  0.603770)  (-2.353923,  0.626921)  (-2.333890,  0.650234)  (-2.313856,  0.673621)
		(-2.293823,  0.696979)  (-2.273790,  0.720192)  (-2.253756,  0.743127)  (-2.233723,  0.765637)  (-2.213689,  0.787555)
		(-2.193656,  0.808696)  (-2.173623,  0.828855)  (-2.153589,  0.847807)  (-2.133556,  0.865303)  (-2.113523,  0.881072)
		(-2.093489,  0.894823)  (-2.073456,  0.906239)  (-2.053422,  0.914983)  (-2.033389,  0.920695)  (-2.013356,  0.922998)
		(-1.993322,  0.921498)  (-1.973289,  0.915789)  (-1.953255,  0.905457)  (-1.933222,  0.890085)  (-1.913189,  0.869265)
		(-1.893155,  0.842605)  (-1.873122,  0.809738)  (-1.853088,  0.770339)  (-1.833055,  0.724141)  (-1.813022,  0.670947)
		(-1.792988,  0.610655)  (-1.772955,  0.543277)  (-1.752922,  0.468963)  (-1.732888,  0.388023)  (-1.712855,  0.300954)
		(-1.692821,  0.208462)  (-1.672788,  0.111485)  (-1.652755,  0.011213)  (-1.632721, -0.090902)  (-1.612688, -0.193133)
		(-1.592654, -0.293486)  (-1.572621, -0.389713)  (-1.552588, -0.479342)  (-1.532554, -0.559723)  (-1.512521, -0.628103)
		(-1.492487, -0.681707)  (-1.472454, -0.717860)  (-1.452421, -0.734118)  (-1.432387, -0.728433)  (-1.412354, -0.699322)
		(-1.392321, -0.646060)  (-1.372287, -0.568858)  (-1.352254, -0.469044)  (-1.332220, -0.349197)  (-1.312187, -0.213244)
		(-1.292154, -0.066474)  (-1.272120,  0.084536)  (-1.252087,  0.232119)  (-1.232053,  0.367838)  (-1.212020,  0.482962)
		(-1.191987,  0.569059)  (-1.171953,  0.618722)  (-1.151920,  0.626355)  (-1.131886,  0.588969)  (-1.111853,  0.506910)
		(-1.091820,  0.384387)  (-1.071786,  0.229724)  (-1.051753,  0.055185)  (-1.031720, -0.123686)  (-1.011686, -0.289270)
		(-0.991653, -0.423599)  (-0.971619, -0.510489)  (-0.951586, -0.537872)  (-0.931553, -0.500019)  (-0.911519, -0.399266)
		(-0.891486, -0.246836)  (-0.871452, -0.062395)  (-0.851419,  0.127893)  (-0.831386,  0.294810)  (-0.811352,  0.410707)
		(-0.791319,  0.454673)  (-0.771285,  0.417286)  (-0.751252,  0.303844)  (-0.731219,  0.135005)  (-0.711185, -0.055846)
		(-0.691152, -0.228579)  (-0.671119, -0.344838)  (-0.651085, -0.377417)  (-0.631052, -0.318133)  (-0.611018, -0.181747)
		(-0.590985, -0.004045)  (-0.570952,  0.166393)  (-0.550918,  0.281304)  (-0.530885,  0.307388)  (-0.510851,  0.237876)
		(-0.490818,  0.096578)  (-0.470785, -0.068212)  (-0.450751, -0.199208)  (-0.430718, -0.250499)  (-0.410684, -0.205472)
		(-0.390651, -0.084676)  (-0.370618,  0.060926)  (-0.350584,  0.170183)  (-0.330551,  0.197940)  (-0.310518,  0.135759)
		(-0.290484,  0.016700)  (-0.270451, -0.099729)  (-0.250417, -0.156642)  (-0.230384, -0.128877)  (-0.210351, -0.036434)
		(-0.190317,  0.066389)  (-0.170284,  0.121600)  (-0.150250,  0.101080)  (-0.130217,  0.022939)  (-0.110184, -0.060426)
		(-0.090150, -0.096204)  (-0.070117, -0.065499)  (-0.050083,  0.005110)  (-0.030050,  0.063474)  (-0.010017,  0.069741)
		( 0.010017,  0.024720)  ( 0.030050, -0.032727)  ( 0.050083, -0.057437)  ( 0.070117, -0.033741)  ( 0.090150,  0.013753)
		( 0.110184,  0.043122)  ( 0.130217,  0.031782)  ( 0.150250, -0.005395)  ( 0.170284, -0.032195)  ( 0.190317, -0.025671)
		( 0.210351,  0.003688)  ( 0.230384,  0.024915)  ( 0.250417,  0.018539)  ( 0.270451, -0.005216)  ( 0.290484, -0.019775)
		( 0.310518, -0.011417)  ( 0.330551,  0.007613)  ( 0.350584,  0.015178)  ( 0.370618,  0.004720)  ( 0.390651, -0.009215)
		( 0.410684, -0.010143)  ( 0.430718,  0.000929)  ( 0.450751,  0.008905)  ( 0.470785,  0.004657)  ( 0.490818, -0.004549)
		( 0.510851, -0.006352)  ( 0.530885,  0.000278)  ( 0.550918,  0.005301)  ( 0.570952,  0.002388)  ( 0.590985, -0.003166)
		( 0.611018, -0.003347)  ( 0.631052,  0.001129)  ( 0.651085,  0.003150)  ( 0.671119,  0.000295)  ( 0.691152, -0.002423)
		( 0.711185, -0.001056)  ( 0.731219,  0.001605)  ( 0.751252,  0.001319)  ( 0.771285, -0.000918)  ( 0.791319, -0.001285)
		( 0.811352,  0.000430)  ( 0.831386,  0.001114)  ( 0.851419, -0.000126)  ( 0.871452, -0.000906)  ( 0.891486, -0.000041)
		( 0.911519,  0.000712)  ( 0.931553,  0.000116)  ( 0.951586, -0.000552)  ( 0.971619, -0.000136)  ( 0.991653,  0.000427)
		( 1.011686,  0.000127)  ( 1.031720, -0.000333)  ( 1.051753, -0.000103)  ( 1.071786,  0.000263)  ( 1.091820,  0.000075)
		( 1.111853, -0.000209)  ( 1.131886, -0.000048)  ( 1.151920,  0.000168)  ( 1.171953,  0.000023)  ( 1.191987, -0.000134)
		( 1.212020, -0.000003)  ( 1.232053,  0.000105)  ( 1.252087, -0.000013)  ( 1.272120, -0.000080)  ( 1.292154,  0.000024)
		( 1.312187,  0.000059)  ( 1.332220, -0.000030)  ( 1.352254, -0.000040)  ( 1.372287,  0.000032)  ( 1.392321,  0.000024)
		( 1.412354, -0.000031)  ( 1.432387, -0.000011)  ( 1.452421,  0.000027)  ( 1.472454,  0.000001)  ( 1.492487, -0.000021)
		( 1.512521,  0.000006)  ( 1.532554,  0.000015)  ( 1.552588, -0.000009)  ( 1.572621, -0.000009)  ( 1.592654,  0.000010)
		( 1.612688,  0.000003)  ( 1.632721, -0.000009)  ( 1.652755,  0.000001)  ( 1.672788,  0.000007)  ( 1.692821, -0.000003)
		( 1.712855, -0.000004)  ( 1.732888,  0.000004)  ( 1.752922,  0.000002)  ( 1.772955, -0.000004)  ( 1.792988,  0.000000)
		( 1.813022,  0.000003)  ( 1.833055, -0.000001)  ( 1.853088, -0.000002)  ( 1.873122,  0.000002)  ( 1.893155,  0.000000)
		( 1.913189, -0.000002)  ( 1.933222,  0.000000)  ( 1.953255,  0.000001)  ( 1.973289, -0.000001)  ( 1.993322, -0.000001)
		( 2.013356,  0.000001)  ( 2.033389, -0.000000)  ( 2.053422, -0.000001)  ( 2.073456,  0.000000)  ( 2.093489,  0.000000)
		( 5.919866, -0.000000)  ( 5.939900, -0.000000)  ( 5.959933,  0.000000)  ( 5.979967,  0.000000)  ( 6.000000, -0.000000)
	};
	\end{axis}
	\end{tikzpicture}
	\caption{The integrand of $\mathcal{I}(s) =  \int_{-\infty}^{\infty} f\left(M\phi(t)\right) \cdot M\phi'(t) \mathrm{d} t$ where $s=0.01$, $\nu=\sfrac{5}{2}$, $n_\gamma=5$, $n_x=0$, $\lambda=0$, $R_1=6.31$, $\zeta_1=1$, $R_2=2.0$, and $\zeta_2=1$ (see row 2 of Table \ref{table:with_phi12}).}
	\label{Transformed_Integrand}
\end{figure}
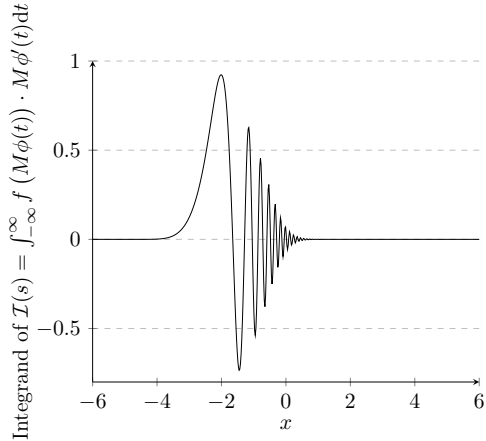

In particular, it is worth noting that the number of collocation points needed to approximate the portion of the integral below zero is often higher than the number of collocation points needed to approximate the portion of the integral above zero. This is due to the asymmetry of the transformed integrand which is represented by the summation in \eqref{approximation}, that is to say of the integrand of $\mathcal{I}(s) =  \int_{-\infty}^{\infty} f\left(M\phi(t)\right) \cdot M\phi'(t) \mathrm{d} t$. This asymmetric behavior can be observed in 	Figure ~\ref{Transformed_Integrand}.

Now in order see a pattern in our data regarding the effect of the collocation points on the accuracy of the approximation, for our first point on the graph, a moderate number of median collocation points was taken in order to ensure sufficient proximity to the exact result. More precisely, for the first point on the graph in Figure~\ref{collocation_points_phi1}, we took the summation in \eqref{approximation} over 70 collocation points, with a lower bound of -62 and an upper bound of 7. From then on, the upper and lower bounds of the summation were each extended by one. We proceeded to extend the bounds of the summation in this manner -- increasing the number of collocation points by two at a time until we reached 138 collocation points, which was the necessary amount of collocation points originally determined by our integrator (see Table~\ref{table:with_phi12}).

For the variable assignments used $\phi_2(t)$ \eqref{phi_2}, our present integrator found $M \approx 21.70135 $ to be an optimal value for $M$; we therefore assigned said value to $M$ to collect our data points. Furthermore, while using transformation $\phi_2(t)$ \eqref{phi_2} and the variable assignments in Figure~\ref{collocation_points_phi1} to approximate $\mathcal{I}(s)$, our present integrator truncated the summation in \eqref{approximation} at -63 and at 33. In a similar way as was done for transformation $\phi_1(t)$, for our first point of data, we took the summation over 45 collocation points, with a lower bound of -37 and an upper bound of 7. From then on, the upper and lower bounds of the summation were each extended by one until we reached a total of 97 collocation points, which was the number of collocation points evaluated by our integrator.

Notice that the graphs in Figure~\ref{collocation_points_phi1} both rapidly decrease and approach zero as the number of collocation points increase. Whether our integrator uses transformation $\phi_1(t)$ \eqref{phi_1} or transformation $\phi_2(t)$ \eqref{phi_2}, the behaviour of the absolute error for the approximation of $\mathcal{I}(s)$ based on the number of collocation points is strikingly similar. In fact, after a certain number of collocation points, both graphs follow a pattern of double exponential decay.

With the particular variable assignments used in Figure~\ref{collocation_points_phi1}, it is clear that using transformation $\phi_2(t)$ \eqref{phi_2}, less collocation points are needed for the absolute error to converge to zero. This is representative of all possible variable assignments for the evaluation of $\mathcal{I}(s)$; with the variable assignments we selected, there are always more collocation points evaluated when using transformation $\phi_1(t)$ \eqref{phi_1} (see Table~\ref{table:with_phi12} ).

That being said however, the parameter $M$ is perhaps just as significant in determining the efficiency of our integrator. The consequences of $M$ on our integrator's efficiency will be addressed in the following subsection.

\subsection{Effect of $M$}
\begin{figure}[H]
	\pgfplotsset{width=16cm, height=8cm}
	\centering
	\begin{tikzpicture}[scale=0.75]
	\begin{axis}[
	axis lines = left,
	xlabel={$M$},
	ylabel={Absolute error ($10^{-5}$)},
	xmin=4, xmax=19,
	ymin=0, ymax=12.5,
	ymajorgrids=true,
	grid style=dashed,
	]
	\addplot[black,thick]
	coordinates {
	(10.0000	,	13.900000000000000)
	(10.1667	,	11.956722253380000)
	(10.3333	,	10.286919976090000)
	(10.5000	,	8.7689369761600000)
	(10.6667	,	7.4023778378400000)
	(10.8333	,	6.1833575998400000)
	(11.0000	,	5.1054436206200000)
	(11.1667	,	4.1604288943600000)
	(11.3333	,	3.3389570499900000)
	(11.5000	,	2.6310201417700000)
	(11.6667	,	2.0263525057200000)
	(11.8333	,	1.5147294603700000)
	(12.0000	,	1.0861980382400000)
	(12.1667	,	0.7312399180000000)
	(12.3333	,	0.4408870991500000)
	(12.5000	,	0.2067976321700000)
	(12.6667	,	0.0212962255700000)
	(12.8333	,	0.1226070925000000)
	(13.0000	,	0.2312203035600000)
	(13.1667	,	0.3101867255100000)
	(13.3333	,	0.3645130047300000)
	(13.5000	,	0.3986055520100000)
	(13.6667	,	0.4163097066900000)
	(13.8333	,	0.4209551117100000)
	(14.0000	,	0.4153996768400000)
	(14.1667	,	0.4020750433200000)
	(14.3333	,	0.3830297750600000)
	(14.5000	,	0.3599731631700000)
	(14.6667	,	0.3343134452200000)
	(14.8333	,	0.3071968278000000)
	(15.0000	,	0.2795407737000000)
	(15.1667	,	0.2520669992600000)
	(15.3333	,	0.2253297950900000)
	(15.5000	,	0.1997419752300000)
	(15.6667	,	0.1755976388700000)
	(15.8333	,	0.1530926924300000)
	(16.0000	,	0.1323442373300000)
	(16.1667	,	0.1134036105100000)
	(16.3333	,	0.0962723487600000)
	(16.5000	,	0.0809116225500000)
	(16.6667	,	0.0672537994300000)
	(16.8333	,	0.0552093995900000)
	(17.0000	,	0.0446740659500000)
	(17.1667	,	0.0355352148900000)
	(17.3333	,	0.0276753795600000)
	(17.5000	,	0.0209753204300000)
	(17.6667	,	0.0153186892400000)
	(17.8333	,	0.0105927643500000)
	(18.0000	,	0.0066895794500000)
	(18.1667	,	0.0000000000000000)
	(19.0000	,	0.0000000000000000)
	};
	\addplot[black,ultra thick,dashed]
	coordinates {
	(4.0000		,	12.463989799841600)
	(4.1667		,	7.7757229544406900)
	(4.3333		,	4.5491581368605100)
	(4.5000		,	2.4615009949393000)
	(4.6667		,	1.1720618816384600)
	(4.8333		,	0.4064313871690200)
	(5.0000		,	0.0287978410720500)
	(5.1667		,	0.2613450731975300)
	(5.3333		,	0.3710286416950800)
	(5.5000		,	0.4048283336320200)
	(5.6667		,	0.3902869913687600)
	(5.8333		,	0.3455580198948600)
	(6.0000		,	0.2845421883385000)
	(6.1667		,	0.2183537209075900)
	(6.3333		,	0.1553744801215600)
	(6.5000		,	0.1011632908284700)
	(6.6667		,	0.0585438500004900)
	(6.8333		,	0.0279676909853700)
	(7.0000		,	0.0081577171094500)
	(7.1667		,	0.0031140772444400)
	(7.3333		,	0.0083313459623500)
	(7.5000		,	0.0097445478809800)
	(7.6667		,	0.0091104331295400)
	(7.8333		,	0.0076231745995100)
	(8.0000		,	0.0059803684521400)
	(8.1667		,	0.0045112217227200)
	(8.3333		,	0.0033177004424000)
	(8.5000		,	0.0023879334519900)
	(8.6667		,	0.0016699673601800)
	(8.8333		,	0.0011147489405300)
	(9.0000		,	0.0006864156274100)
	(9.1667		,	0.0003640877169500)
	(9.3333		,	0.0001336245843800)
	(9.5000		,	0.0000177898987700)
	(9.6667		,	0.0000000000000000)
	(19.0000	,	0.0000000000000000)
	};
\legend{$\phi_1(t)$, $\phi_2(t)$}
	\end{axis}
	\end{tikzpicture}
	\caption{Absolute error dependence of $\mathcal{I}_M(s)$ \eqref{approximation} on the total number of collocation points using transformations $\phi_1(t)$ \eqref{phi_1} and $\phi_2(t)$ \eqref{phi_2}, and the variable assignments from the second row of Table~\ref{table:with_phi12}. The absolute error was found by comparing the results of the approximation with the result obtained by a MATLAB built-in numerical integration function that uses global adaptive quadrature.}
\label{figure:M_error}
\end{figure}
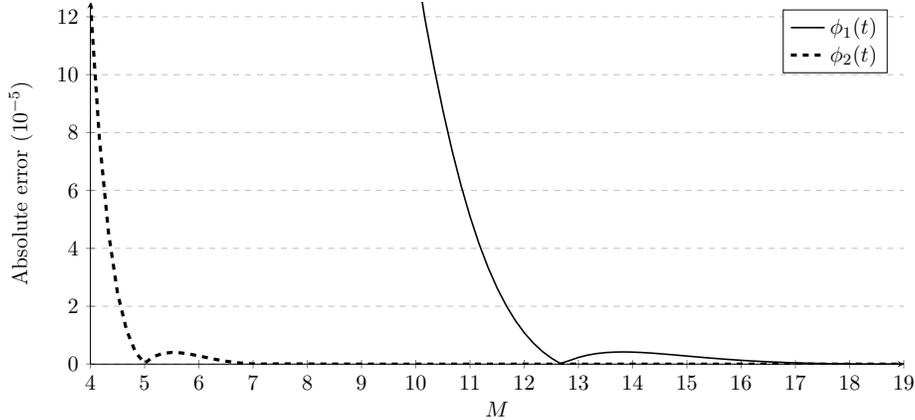

In Figure~\ref{figure:M_error}, it can easily be seen that the absolute error converges to zero much more quickly using $\phi_2(t)$ \eqref{phi_2}. In other words, for a given set of variable assignments, a smaller mesh size is required using transformation $\phi_1(t)$ \eqref{phi_1} to obtain results that are as accurate as the results using $\phi_2(t)$ \eqref{phi_2}.

The challenge is to find a value for $M$ that is large enough such that the error converges to zero, yet small enough so as to not make $h$, the mesh size, prohibitively small. Clearly, given \eqref{mh=lambda}, as $M$ increases in size, $h$ decreases. For this reason, it would be counterproductive to let $M$ be excessively large, as quite a number of collocation points in the summation \eqref{approximation} would be needed to compensate for the unnecessarily small mesh size. Yet if $M$ is too small, as shown in Figure~\ref{figure:M_error}, the error of the approximation will not approach zero. This is what is meant by finding an ``optimal value'' for $M$.

Our present integrator finds an optimal value for $M$ by first letting $M$ be a relatively small positive constant, say $M_1$, based on the relative error tolerance input by the user and according to the following formula~\cite{Ooura-Mori-112-229-99}:
\begin{equation}
	M_1 \,=\, \frac{-\pi}{A} \, \log\left(\sqrt{\varepsilon_0}\right),
\end{equation}
and $\varepsilon_0$ denotes the relative error tolerance. For the transformation $\phi_1(t)$ \eqref{phi_1}, we use $A = 2$, and for the transformation $\phi_2(t)$ \eqref{phi_2}, we use $A = 5$.  For details regarding the selection of each subsequent $M_i$, we direct the reader to the algorithm in \cite{Ooura-Mori-112-229-99}.

The integrator then calculates $\mathcal{I}_{M_1}(s)$ using equation \eqref{approximation}, with the value of $M_1$ assigned to $M$. The result is not a viable approximation in the sense that it is quite far away from the desired result due to the fact that the value assigned to $M$ is too small to expect the error to converge to zero. Rather, $\mathcal{I}_{M_1}(s)$ serves as a value to which a second, more accurate approximation can be compared for the purpose of determining a relative error.

For the second evaluation of \eqref{approximation}, we assign a larger and more suitable value to $M$, say $M_2$. The ensuing approximation, $\mathcal{I}_{M_2}(s)$, often times meets the user's relative error tolerance requirements, at which point the integrator outputs the value of $\mathcal{I}_{M_2}(s)$ as its final result. Otherwise, the process must be repeated with another value assigned to $M$, until a result of sufficient accuracy is found. This could take up to a maximum of four different assignments to $M$, each of which is larger than the previous.

Clearly, with this type of integrator, it is ideal to find an optimal value for $M$ on the second attempt, and highly undesirable to not find an optimal value until the fourth attempt, because the approximation \eqref{approximation} would need to be calculated four times over, meaning that numerous time-consuming computations would have to be executed. It is for this reason that $M$ also holds significant weight while comparing the efficiency of our algorithm using transformations $\phi_1(t)$~\eqref{phi_1} and $\phi_2(t)$~\eqref{phi_2}.
For instance, by solely comparing the collocation points in Table~\ref{table:with_phi12}, it may appear as though $\phi_2(t)$ \eqref{phi_2} is without fail the more efficient transformation to use. That being said, since the value $n_M$, which we are using to denote the total number of values that are assigned to $M$ before successfully completing the approximation for $\mathcal{I}(s)$, also gives valuable insight into the efficiency of each transformation. A high $n_M$ value  signifies that our integrator had to complete the approximation \eqref{approximation} several times before finding an optimal value for $M$ and successfully approximating the integral, thus multiplying the total number of calculations that needed to be computed by the integrator.

Observing the sixth row in Table~\ref{table:with_phi12}, for example, it can be seen that although the summation in \eqref{approximation} is truncated with 81 collocation points using transformation $\phi_1(t)$ \eqref{phi_1}, and only 72 collocation points using transformation $\phi_2(t)$ \eqref{phi_2}, the values for $n_M$ are found to be 2 and 4, using transformation $\phi_1(t)$ \eqref{phi_1} and $\phi_2(t)$ \eqref{phi_2}, respectively. In this particular case, therefore, it can be argued that our present integrator is more efficient in its approximation using transformation $\phi_1(t)$ \eqref{phi_1}. Generally speaking, however, the two-parameter transformation, $\phi_2(t)$ \eqref{phi_2} lends itself to a more efficient approximation of the semi-infinite spherical Bessel integrals.

\section{Conclusion}
The presence of the semi-infinite spherical Bessel integral in the analytical expression for three-centre nuclear attraction integrals causes the rapid and accurate numerical evaluation of the integrals to be extremely challenging. The $S$ transformation, introduced by Safouhi, is able to transform the spherical Bessel integrals into considerably simpler sine integrals. The zeros of the integrand become equidistant, which significantly increases our ability to apply various extrapolation and transformation methods.  Nevertheless, due to its slow convergence, the semi-infinite sine integrals still pose substantial computational difficulties.

In the present work, we presented an efficient way of treating the arduous sine integrals that result from the $S$ transformation. Utilizing the slowly decaying, highly oscillatory nature of the integrands, we applied two distinct double exponential transformations. These transformations each rendered highly efficient quadrature formulae that spanned from negative to positive infinity. That being said, we were able to truncate the summations at relatively small positive and negative integers. It was found that remarkably few collocation points were required to obtain sufficiently accurate results. This was especially the case for the second transformation. In our numerical section, it can be seen that by using the second transformation, the semi-infinite spherical Bessel integral can be approximated with a high predetermined accuracy in under 100 collocation points. The numerical tests served to further verify the precision and demonstrate the significance of the method.

\section{Numerical tables}

\clearpage

\begin{sidewaystable}[H]
\begin{center}
\setlength{\tabcolsep}{3pt}
\caption{Computation of $\mathcal{I}(s)$. \label{table:with_phi12}}
\begin{tabular*}{\hsize}{@{\extracolsep{\fill}}ccccccccccccccccccc}			\hline			\\[-.9em]
$s$ & $\nu$ & $n_\gamma$ & $n_x$ & $\lambda$ & $R_1$ & $\zeta_1$ & $R_2$ & $\zeta_2$ &$\mathcal{I}(s)^{\textrm{Matlab}}$ & $n^{\phi_1}$ & $\max^{\phi_1}$ & $n_M^{\phi_1}$ & $\varepsilon_M^{\phi_1}$
& $n^{\phi_2}$ & $\max^{\phi_2}$ & $n_M^{\phi_2}$  & $\varepsilon_M^{\phi_2}$\\
\hline \\[-.9em]
0.99 & $\sfrac{5}{2}$  & 1 & 0 & 0 & 24.00  & 1.5  & 2.0 & 1 & .113 874 170 637 205(~0) & 142 & 45 & 2 & .1(-16) & 93 & 34 & 2 & .1(-16)   \\ [+.25em]
0.01 & $\sfrac{5}{2}$  & 5 & 0 & 0 & ~6.31  & 1.0  & 2.0 & 1 & .638 243 453 884 443(~0) & 138 & 41 & 2 & .1(-16) & 97 & 33 & 2 & .1(-16)   \\ [+.25em]
0.99 & $\sfrac{5}{2}$  & 5 & 0 & 0 & ~4.50  & 2.0  & 1.5 & 1 & .701 581 269 512 308(~0) & 139 & 42 & 2 & .1(-16) & 97 & 33 & 2 & .1(-16)   \\ [+.25em]
0.99 & $\sfrac{5}{2}$  & 5 & 1 & 0 & ~3.00  & 1.5  & 3.5 & 2 & .242 778 918 544 382(-3) & ~78 & 41 & 2 & .1(-16) & 78 & 37 & 4 & .1(-16)   \\ [+.25em]
0.99 & $\sfrac{9}{2}$  & 9 & 1 & 1 & ~6.00  & 2.0  & 3.5 & 1 & .183 138 910 224 195(~1) & 139 & 42 & 2 & .1(-16) & 97 & 33 & 2 & .1(-16)   \\ [+.25em]
0.01 & $\sfrac{9}{2}$  & 9 & 2 & 1 & ~8.50  & 2.0  & 3.5 & 2 & .248 336 723 989 982(-3) & ~81 & 44 & 2 & .1(-16) & 72 & 35 & 4 & .5(-18)   \\ [+.25em]
0.01 & $\sfrac{7}{2}$  & 3 & 2 & 1 & ~3.00  & 2.0  & 5.0 & 1 & .285 091 100 421 789(-2) & ~85 & 41 & 3 & .1(-16) & 81 & 37 & 3 & .1(-16)   \\ [+.25em]
0.99 & $\sfrac{7}{2}$  & 5 & 2 & 2 & ~4.00  & 2.5  & 5.5 & 1 & .112 567 767 257 153(~0) & 139 & 42 & 2 & .1(-16) & 92 & 33 & 2 & .1(-16)   \\ [+.25em]
0.01 & $\sfrac{9}{2}$  & 7 & 2 & 2 & ~9.00  & 2.0  & 3.5 & 1 & .183 269 571 025 263(-2) & 141 & 44 & 2 & .1(-16) & 97 & 33 & 2 & .1(-16)   \\ [+.25em]
0.01 & $\sfrac{13}{2}$ & 9 & 3 & 2 & ~4.00  & 2.5  & 5.5 & 1 & .167 566 737 865 368(-1) & ~91 & 37 & 3 & .1(-16) & 92 & 36 & 3 & .1(-16)   \\
\hline
\end{tabular*}
\begin{tablenotes}
\begin{itemize}
	\item $\mathcal{I}(s)^{\textrm{Matlab}}$ are obtained by means of a MATLAB built-in numerical integration function that uses global adaptive quadrature.
	\item $\varepsilon_M^{\phi_1}$ stand for the relative errors obtained by using the approximation \eqref{approximation} with the transformation $\phi_1(t)$ \eqref{phi_1}.
	\item $\varepsilon_M^{\phi_2}$ stand for the relative errors obtained by using the approximation \eqref{approximation} with the transformation $\phi_2(t)$ \eqref{phi_2}.
	\item $n^{\phi_1}$ and $n^{\phi_2}$ represent the number of collocation points needed to approximate the integral using \eqref{approximation} with transformations $\phi_1(t)$ \eqref{phi_1} and $\phi_2(t)$ \eqref{phi_2}, respectively.
	\item $n_M^{\phi_1}$ and $n_M^{\phi_2}$ represent the total number of values that are assigned to $M$ before successfully completing the approximation \eqref{approximation} with transformations $\phi_1(t)$ \eqref{phi_1} and $\phi_2(t)$ \eqref{phi_2}, respectively.
	\item $\max^{\phi_1}$ and $\max^{\phi_2}$ represent the upper limit of the summation in \eqref{approximation} with transformations $\phi_1(t)$ \eqref{phi_1} and $\phi_2(t)$ \eqref{phi_2}, respectively.
	\item The approximation \eqref{approximation} was implemented in Python with symbolic computation package, SymPy
	\item The error tolerance of $\varepsilon = 10^{-15}$ was used in our calculation.
	\item Calculations were completed using IEEE 754 double precision
 	\item Numbers in parentheses represent powers of 10.
\end{itemize}
\end{tablenotes}
\end{center}
\end{sidewaystable}

\begin{sidewaystable}[H]
\begin{center}
\setlength{\tabcolsep}{3pt}
\caption{Computation of $\mathcal{I}(s)$. \label{table:using_phi12}}
\begin{tabular*}{\hsize}{@{\extracolsep{\fill}}ccccccccccccccccccccccccc}			\hline			\\[-.9em]
$s$ & $\nu$ & $n_\gamma$ & $n_x$ & $\lambda$ & $R_1$ & $\zeta_1$ & $R_2$ & $\zeta_2$ &$\mathcal{I}(s)^{\phi_1}$ & $n^{\phi_1}$ & $\max^{\phi_1}$ & $n_M^{\phi_1}$ & $\varepsilon_M^{\phi_1}$& $\mathcal{I}(s)^{\phi_2}$ & $n^{\phi_2}$ & $\max^{\phi_2}$ & $n_M^{\phi_2}$  & $\varepsilon_M^{\phi_2}$\\
\hline \\[-.9em]
0.99 &	$\sfrac{9}{2}$  & ~9  & 2 & 1 & 35 & 2.5 & 3.5 & 0.5 & .737 829 982 455 112(~0) & ~82 & 45 & 2 & .1(-16) & .737 829 982 455 148(~0) & 84 & 41 & 4 & .1(-16) \\ [+.25em]
0.01 &	$\sfrac{9}{2}$  & ~3  & 2 & 2 & 45 & 1.5 & 2.0 & 1.0 & .103 851 188 766 556(-2) & 142 & 45 & 2 & .1(-16) & .103 851 188 766 626(-2) & 93 & 34 & 2 & .1(-16) \\ [+.25em]
0.99 &	$\sfrac{13}{2}$ & ~5  & 3 & 3 & 50 & 1.5 & 2.0 & 1.0 & .317 648 310 603 089(-1) & 141 & 44 & 2 & .1(-16) & .317 648 310 603 304(-1) & 93 & 34 & 2 & .1(-16) \\ [+.25em]
0.01 &	$\sfrac{15}{2}$ & ~6  & 4 & 3 & 55 & 1.5 & 2.0 & 1.0 & .989 261 101 060 361(-3) & ~83 & 46 & 2 & .1(-16) & .989 261 101 060 357(-3) & 86 & 42 & 4 & .5(-20) \\ [+.25em]
0.99 &	$\sfrac{15}{2}$ & ~6  & 4 & 4 & 55 & 1.5 & 2.0 & 1.0 & .366 551 897 499 993(-1) & 141 & 44 & 2 & .1(-16) & .366 551 897 500 241(-1) & 93 & 34 & 2 & .5(-16) \\ [+.25em]
0.01 &	$\sfrac{17}{2}$ & ~9  & 5 & 4 & 55 & 1.5 & 2.0 & 1.0 & .698 018 626 122 216(-3) & ~83 & 46 & 2 & .1(-16) & .698 018 626 122 213(-3) & 82 & 40 & 4 & .2(-19) \\ [+.25em]
0.99 &	$\sfrac{23}{2}$ & 21  & 6 & 5 & 55 & 1.5 & 2.0 & 1.5 & .564 022 921 688 577(-2) & ~81 & 44 & 2 & .1(-16) & .564 022 921 688 556(-2) & 84 & 41 & 4 & .1(-16) \\ [+.25em]
0.01 &	$\sfrac{27}{2}$ & 29  & 6 & 6 & 60 & 1.5 & 2.0 & 1.0 & .414 131 181 249 046(~0) & 142 & 45 & 2 & .1(-16) & .414 131 181 249 327(~0) & 93 & 34 & 2 & .1(-16) \\ [+.25em]
0.01 &	$\sfrac{29}{2}$ & 25  & 7 & 6 & 55 & 2.0 & 3.0 & 2.0 & .284 952 750 728 949(-2) & ~82 & 45 & 2 & .1(-16) & .284 952 750 728 952(-2) & 82 & 40 & 4 & .1(-16) \\ [+.25em]
0.01 &	$\sfrac{31}{2}$ & 23  & 7 & 7 & 55 & 2.5 & 2.0 & 1.0 & .107 097 615 409 941(-1) & 143 & 45 & 2 & .1(-16) & .107 097 615 410 016(-1) & 93 & 34 & 2 & .1(-16) \\ [+.25em]
0.01 &	$\sfrac{33}{2}$ & 33  & 7 & 7 & 65 & 2.0 & 2.0 & 1.0 & .167 421 970 712 946(-1) & 143 & 45 & 2 & .1(-16) & .167 421 970 713 064(-1) & 93 & 34 & 2 & .1(-16) \\
\hline
\end{tabular*}
\begin{tablenotes}
\begin{itemize}
	\item $\mathcal{I}(s)^{\textrm{Matlab}}$ are obtained by means of a MATLAB built-in numerical integration function that uses global adaptive quadrature.
	\item $\varepsilon_M^{\phi_1}$ stand for the relative errors obtained by using the approximation \eqref{approximation} with the transformation $\phi_1(t)$ \eqref{phi_1}.
	\item $\varepsilon_M^{\phi_2}$ stand for the relative errors obtained by using the approximation \eqref{approximation} with the transformation $\phi_2(t)$ \eqref{phi_2}.
	\item $n^{\phi_1}$ and $n^{\phi_2}$ represent the number of collocation points needed to approximate the integral using \eqref{approximation} with transformations $\phi_1(t)$ \eqref{phi_1} and $\phi_2(t)$ \eqref{phi_2}, respectively.
	\item $n_M^{\phi_1}$ and $n_M^{\phi_2}$ represent the total number of values that are assigned to $M$ before successfully completing the approximation \eqref{approximation} with transformations $\phi_1(t)$ \eqref{phi_1} and $\phi_2(t)$ \eqref{phi_2}, respectively.
	\item $\max^{\phi_1}$ and $\max^{\phi_2}$ represent the upper limit of the summation in \eqref{approximation} with transformations $\phi_1(t)$ \eqref{phi_1} and $\phi_2(t)$ \eqref{phi_2}, respectively.
	\item The approximation \eqref{approximation} was implemented in Python with symbolic computation package, SymPy
	\item The error tolerance of $\varepsilon = 10^{-15}$ was used in our calculation.
	\item Calculations were completed using IEEE 754 double precision
 	\item Numbers in parentheses represent powers of 10.
\end{itemize}
\end{tablenotes}
\end{center}
\end{sidewaystable}

\clearpage

\end{document}